\documentclass[a4,12pt]{article}
\usepackage{authblk}
\usepackage{amscd}
\usepackage{amssymb}
\usepackage{amsmath}
\usepackage{latexsym}
\def\proof{\noindent{\bf Proof.}\quad}

\usepackage{amsmath}
\def\endproof{\quad {\bf q.e.d.}\\\smallskip}

\newtheorem{theorem}{Theorem}

\newtheorem{proposition}{Proposition}

\title{Principal Moment Analysis}

\date{}
\author[1,2,3]{Magnus Fontes\thanks{fontes.magnus@gene.com; Corresponding author}}
\author[4]{Rasmus Henningsson}

\affil[1]{{\small Department of Cancer Immunology, Genentech, South San Francisco, USA}}
\affil[2]{{\small The Center for Genomic Medicine, Rigshospitalet, Copenhagen, Denmark}}
\affil[3]{{\small Persimune, The Centre of Excellence for Personalized Medicine, Copenhagen, Denmark}}
\affil[4]{{\small Division of Clinical Genetics, Lund University, Lund, Sweden}}

\begin{document}
\maketitle
\section*{Abstract}
Principal Moment Analysis is a method designed for dimension reduction, analysis and visualization of high dimensional multivariate data.
It generalizes Principal Component Analysis and allows for significant statistical modeling flexibility, when approximating an unknown underlying probability distribution, by enabling direct analysis of general approximate measures.
Through  https://principalmomentanalysis.github.io/ we provide an implementation, together with a graphical user interface, of a simplex based version of Principal Moment Analysis.
\section*{Introduction}
Principal Moment Analysis (PMA) is designed for dimension reduction, analysis and visualization of high dimensional multivariate data.
PMA generalizes Principal Component Analysis (PCA) introduced by Pearson and Hotellinger in the beginning of the 20th century (\cite{pea}, \cite{ho1}). 

We start from the observation that, given a set of sample points in feature space, PCA can be formulated as a spectral decomposition of an operator represented by a second moment tensor of a measure defined by a sum of Dirac delta measures corresponding to the given sample points. The usefulness of PCA stems from the fact that this spectral decomposition enables exploration of, in a certain sense, optimal lower dimensional representations of the data. 
This particular formulation of PCA suggests a natural generalization and PMA generalizes PCA  
by allowing us to replace the sum of Dirac delta measures of PCA by general measures that, through some specified design, have been constructed from the given data. This measure construction step in the  PMA analysis scheme permits substantial statistical modeling flexibility and it is followed by a spectral decomposition of a second moment tensor connected with the measure. This spectral decomposition enables exploration of optimal lower dimensional representations of both the measure and the original data.

Through https://principalmomentanalysis.github.io/ we provide an implementation of a simplex based version of PMA, {\it simplex PMA} where we, given a set of sample points in feature space, construct measures based on sums of Hausdorff measures with support on simplexes. The simplexes in a sum can be of varying dimension and are spanned by specified convex combinations of the underlying sample points in feature space. The specification of which samples that contribute to the Hausdorff measure of a simplex can in a flexible way be assigned using the given data and metadata annotations.
The assignation can e.g. be based on that the samples are nearest neighbors, are close together, or be based on that they share the same value for a given metadata annotation. As an illustration, in biomedical applications, such an annotation can e.g. represent replicates, biological subtype, multiple biopsies from the same tissue, multiple samples from the same donor or constitute a time series for a particular subject. 

Our implementation of simplex PMA supports {\it R, Julia} and simple text file format based input data frames holding data and accompanying metadata and includes a graphical user interface enabling flexible construction of the underlying simplexes,  visualizations and extraction of simplex PMA features. 
The {\it algorithmic efficiency} of  our simplex PMA implementation is of the same order as corresponding PCA.

\section*{Statistical learning and PMA}
In multivariate statistical data analysis, we think of a given data set as sampling from the {\it law} of a random variable $X$ taking values in a Hilbert space $H$. The law of $X$ is in general given by a {\it push-forward} of an underlying probability measure $P$ by $X$. The push-forward measure $X_*(P) =: \nu $ is in itself a probability measure defined on {\it feature space}, the Hilbert space $H$ where $X$ takes values. 

In statistical inference and decision theory the underlying measure $\nu$ is, at least partly, assumed to be unknown and based on the available data we try to learn as much as possible about it. 
In the PMA analysis framework we interpret this as, based on available samples, trying to construct an informative approximate measure $\mu$ of  $ \nu $. An approximate measure $\mu$ that we then in turn can analyze and visualize using a spectral decomposition of an operator, represented by a second moment tensor of $\mu$, acting on $H$.
We emphasize that, in the case when we first approximate the underlying unknown probability measure $ \nu $ with a normalized finite sum of Dirac point masses $ \mu :=\frac{1}{|I|} \sum_{x_k \in I} \delta_{{\bf x}_k}$ constructed from a given set of sample points $\{x_k\}_{k \in I}$, and then perform a spectral decomposition of the associated second moment tensor operator, the PMA scheme is equivalent to classical Principal Component Analysis (PCA). 
But again, the general PMA analysis framework adds significant additional modeling flexibility by allowing incorporation of general construction schemes to create approximate measures $\mu$ from given data.
In the implementation of {\it simplex PMA}, available through https://principalmomentanalysis.github.io/, we create approximate measures using sums of uniform Hausdorff measures with support on simplexes spanned by specified convex combinations of the given sample points. This can be interpreted as using a kind of simplex based {\it multidimensional  histograms} constructed from the sample point cloud to approximate $\nu$.
Under fairly mild assumptions on the regularity of the underlying probability measure $\nu$ and with almost any reasonable metric on the space of measures on $H$, e.g. the {\it total variation distance} or a distance based on {\it optimal transport}, a normalized sum of general Hausdorff measures with support on multidimensional simplexes can give a better approximation of $\nu$ than  restricting the approximation to normalized sums of sample point masses, i.e. allowing only zero dimensional Hausdorff measures.

Now, given a measure $\mu$ with finite mass and finite second moment in a Hilbert space $H$, the {\it second moment operator of $\mu$} is well defined and is given by $T \circ T^* : H \longrightarrow H$, where  $T^*: H \longrightarrow L^2(H, \mu)$ denotes the inverse {\it Riesz mapping}
$
T^* (x) (\cdot) = (\cdot, x)_H \, ,
$
and  $T: L^2(H, \mu) \longrightarrow H$  is the Hilbert space adjoint of the inverse Riesz mapping.
The spectral decomposition step of the PMA analysis scheme corresponds to spectral decomposition of the second moment operator.
This spectral decomposition in turn provides us with a (dual) Singular Value Decomposition (SVD) of the pair $T$ and $T^*$, that we in the PMA analysis scheme use to identify, analyze and visualize optimal, with respect to second moments, projected measures $\Pi_* (\mu)$ as well as investigate "directions", {\it principal moment axes}, in $H$ that are the principal contributers to these second moments. 

The described SVD scheme works for any triplet of Hilbert space, Borel measure (finite with finite second moment) and corresponding second moment operator, but in the PMA analysis framework the scheme is applied to measures that, through some design, are constructed from given samples.

Often, in applications of PMA, e.g. in simplex PMA, the underlying Hilbert space $H$ is finite dimensional, but the extrinsic  dimension of the feature space $H$ can be very high. On the other hand, the intrinsic dimension of the state space of interest, where the data essentially lives, will often be of much lower dimension, making the dimension reduction step of PMA interesting and natural.
We also point out that the support of the unknown measure $\nu$ and of our approximate measures $\mu$, do not apriori have to be assumed to be well behaved objects, like e.g. locally defined Riemannian manifolds embedded in $H$, for the PMA analysis scheme to make sense. 
The supports of both $\nu$ and $\mu$ can in fact be allowed to be any (measurable) set in $H$, e.g. a set having highly variable local Hausdorff dimension.

Finally, the PMA machinery fits into the framework of {\it Reproducing Hilbert Space Kernels} and the PMA analysis scheme works in the nonlinear Kernel setting, replacing the scalar product $(\cdot,\cdot)_H$ in $H$ by a given nonlinear positive semidefinite Kernel $\kappa (\cdot,\cdot)$, defined on $H \times H$, but this direction will not be pursued here.

\section*{The mathematical framework for PMA}

We begin by establishing some notions and notations, and at the same time highlighting some useful linear algebra and functional analysis results.
Throughout, $H$ will denote a real separable Hilbert space with scalar product $(\cdot , \cdot)_H$. In many applications, $H$ will be finite dimensional, but we present the PMA analysis framework for general Hilbert spaces $H$ since it goes through unchanged and since it will be important when for example generalizing PMA to a nonlinear kernel setting.

A Hilbert space $H$ is by definition a complete metric space with metric defined by a norm $\| \cdot \|_H$ coming from a scalar product
\begin{equation}
\| \cdot \|^2_{H} := (\cdot , \cdot)_H \, .
\end{equation}
We will sometimes leave out the symbol "$H$" in formulas if its clear from the context.
Every separable Hilbert space $H$ has a complete orthonormal countable basis. 
If $\{{\bf e}_k\}_{k=1}^{\infty}$ is such a Hilbert space basis for $H$, then every element ${\bf u} \in H$ can be represented by its sequence of coefficients $\{({\bf u}, {\bf e}_k)\}_{k=1}^\infty$ as
\begin{equation}
{\bf u}= \sum_{k=1}^{\infty} ({\bf u}, {\bf e}_k) {\bf e}_k 
\end{equation}
with convergence in norm, i.e. the mapping $H \ni {\bf u} \mapsto \{({\bf u}, {\bf e}_k)\}_{k=1}^{\infty}$ is a Hilbert space isomorphism from $H$ to $l^2({\bf N})$.

A linear map $ \Pi :H \longrightarrow H$ is an {\it orthogonal projection} if and only if $ \Pi^2 = \Pi$ and $\Pi^* = \Pi$.
Let $\mathcal{P}^s(H)$ denote the set of orthogonal projections on $H$ of rank $s \in {}\bf N$, i.e. orthogonal projections with $s$-dimensional image set.
Note that any $ \Pi \in \mathcal{P}^s(H)$ has finite rank and thus is a compact self-adjoint operator, so by spectral decomposition, given  $ \Pi \in \mathcal{P}^s(H)$, there exists an orthonormal set of vectors ${\bf \pi}_k \in H$, $k=1, \dots , s$ such that
\begin{equation}
\Pi = \sum_{k=1}^s {\bf \pi}_k \otimes {\bf \pi}_k \, .
\end{equation}
This representation is unique up to the orthogonal group acting on the {\it image set} of $\Pi$, equal to the linear hull of  ${\bf \pi}_k \in H$, $k=1, \dots , s$.

As we pointed out, in applications the Hilbert space will often be finite dimensional, i.e. Hilbert space isomorphic to ${\bf R}^p$ for a suitable $p \in {\bf N}$ and we will then identify an element $ {\bf x} \in {\bf R}^p$ with its column vector of coefficients in the natural basis of ${\bf R}^p$, so that we write ${\bf x} = [ x_1, \dots , x_p]^T$ where $T$ denotes matrix transposition and, using matrix multiplication,
$$
({\bf x},{\bf y})_{{\bf R}^p}  =  [ x_1, \dots , x_p] [ y_1, \dots , y_p]^T = : {\bf x}^T {\bf y}\, ,
$$ 
with ${\bf x}, {\bf y} \in {\bf R}^p$.

\subsection*{Measures with finite moments and SVD}

Let $\mu$ be a positive Borel measure on a Hilbert space $H$ with finite zero, first and second order moments, i.e.
\begin{eqnarray}
\int_H \left( 1 + \| {\bf x} \|_H^2 \right) d \mu ({\bf x}) \; < + \infty
\end{eqnarray}
We can then define a bounded linear operator
$$
T: L^2(H, \mu) \mapsto H
$$
by
\begin{equation}\label{eq:definitionT}
(T(u),{\bf v})_H := \int_{H} u({\bf x}) ({\bf x},{\bf v})_H  d \mu({\bf x}) \, .
\end{equation}
By expansion in an orthonormal basis for $H$, since $\mu$ has finite second moment, it is easy to see that $T$ can be approximated in Hilbert space norm by finite rank operators, i.e. $T$ is a compact operator.

From (\ref{eq:definitionT}) we can conclude that the Hilbert space adjoint operator $T^*: H \mapsto L^2(H, \mu)$ is given by

\begin{equation}
T^*({\bf v})({\bf x})= ( {\bf x}, {\bf v})_{H} \, .
\end{equation}

For $ T^* \circ T :  L^2(H, \mu) \mapsto  L^2(H, \mu)$ we get

$$
\left( T^* \circ T (u), v \right)_{L^2(H, \mu)} = \iint_{H \times H} ( {\bf x} , {\bf y} )_H u({\bf x}) v({\bf y})  d \mu({\bf x}) d \mu({\bf y}) \, ,
$$
so
\begin{equation}
T^* \circ T (u) ({\bf y}) = \int_{H} ( {\bf x} , {\bf y} )_H u({\bf x}) d \mu({\bf x}) \, .
\end{equation}
Similarly from
$$
(T^*({\bf u}(\cdot)),T^*({\bf v}(\cdot)  ))_{L^2(H, \mu)} = \int_{H} ({\bf x}, {\bf u})_H ({\bf x}, {\bf v})_H d \mu ({\bf x}) \, ,
$$
we conclude that
\begin{equation}
 (T \circ T^* ({\bf u}), {\bf v})_H = \int_{H} ({\bf x}, {\bf u})_H ({\bf x}, {\bf v})_H d \mu 
 =: (\int_{H} ({\bf x} \otimes {\bf x} ) \,\, d \mu ({\bf x}) \, \llcorner {\bf u} , {\bf v} )_H \, .
\end{equation}

Since $T$ and $T^*$ are compact they have a (dual) Singular Value Decomposition (SVD). 
The spectral decomposition step of PMA analysis scheme builds on this dual SVD of $T$ and $T^*$ and it is achieved through spectral decompositions of either of the compact self-adjoint operators
\begin{equation}
T^* \circ T (u) (\cdot) = \int_{H} ( {\bf x} , \cdot) u({\bf x}) \, d \mu({\bf x}) 
\end{equation}
or
\begin{equation}
T \circ T^* ({\bf u}) = \int_{H} {\bf x} \otimes {\bf x} \, d \mu ({\bf x}) \llcorner {\bf u}  \, .
\end{equation}

 It is natural to, in analogy with classical mechanics, introduce the  {\it moment tensors} of $\mu$. 
 The first moment tensor is defined as
\begin{equation}
{ \bf M}_1 (\mu) := \int_{H} {\bf x} \, d\mu ({\bf x}) 
\end{equation}
 and the second moment tensor of $\mu$
\begin{equation}
{ \bf M_2} (\mu) =  \int_{H} {\bf x} \otimes {\bf x} \, d\mu ({\bf x}) \, .
\end{equation}
Since $\mu$ has finite second moment, the second moment tensor of $\mu$ is a positive trace class operator on $H$ with {\it trace norm} given by
$$
 \| {\bf M}_2((\mu)) \|_{Tr}:  = tr \left( {\bf M}_2((\mu))\right) \, .
$$
We use the trace norm to define optimal projections of given rank of $\mu$.

\begin{proposition}
Given a finite positive Borel measure $\mu$ on $H$, with finite second moment and a fixed rank $s$, there exists a $ \Pi _s \in \mathcal{P}^s(H) $ such that
\begin{equation}\label{eq:optimal}
 \| {\bf M}_2({\Pi_s}_*(\mu)) \|_{Tr} =
\sup_{\Pi  \in \mathcal{P}^s(H)} \| {\bf M}_2(\Pi_*(\mu)) \|_{Tr} \, ,
\end{equation}
where $\| \cdot \|_{Tr}$ denotes the {\it trace norm}.
\end{proposition}
This proposition follows from the spectral decomposition of the compact, self adjoint and positive operator ${\bf M}_2(\mu)$ and solutions, represented by projections $\Pi_s$, constructed based on the $s$ "first" eigenvectors, provides lower dimensional optimal representations of the measure $\mu$.

In the next section we give explicit formulas in the special case when $H={\bf R}^p$, describing the corresponding SVD, the principal moment axes for the operators $T$ and $T^*$ and the corresponding optimal lower dimensional representations of $\mu$.

\subsection*{Finite dimensional feature space}

Recall that
\begin{equation}
(T \circ T^* ({\bf u}), v) = ({ \bf M_2} (\mu), u \otimes v) \, .
\end{equation}
We perform the dual SVD of  $T$ and $T^*$, starting from the symmetric positive semidefinite tensor ${\bf M}_2 (\mu)$, that we, in the finite dimensional case, identify with a symmetric positive semidefinite $p \times p$ matrix. So, with a slight abuse of notation:
\begin{equation}
{\bf M}_2 (\mu) = \int_{{\bf R}^p} {\bf x} {\bf x}^T \, d\mu ({\bf x}) \, .
\end{equation}
The trace norm of ${\bf M}_2 (\mu)$ in this setting can be computed as
\begin{equation}
\| {\bf M}_2(\mu) \|_{Tr}   = \text{trace}  \left( \int_{{\bf R}^p} {\bf x} {\bf x}^T \, d\mu ({\bf x}) \right) = \int_{{\bf R}^p} {\bf x}^T {\bf x} \, d\mu ({\bf x}) \, .
\end{equation}

Now let $\lambda_1 \geq \lambda_2 \geq \cdots \geq \lambda_r > 0 $ be the ordering of the, necessarily positive, non-zero eigenvalues of  ${\bf M}_2(\mu)$. The rank of ${\bf M}_2(\mu)$ is $r$ and the dimension of the null space is $p-r$.

Let ${\bf v}_k$, $k=1, \dots , r$ be corresponding orthogonal eigenvectors, the {\it Principal Moment Axes}, PMAs, i.e.
\begin{equation}
{\bf M}_2 (\mu) {\bf v}_k = \lambda_k {\bf v}_k \; ; \; k=1, \dots , r \, .
\end{equation}

We define dual PMAs, or {\it principal moment functionals}, as follows
$$
u_k({\bf x}):= \lambda_k^{-\frac{1}{2}} ({\bf x}, {\bf v}_k) \quad : \; k=1, \dots , r \, .
$$
We note that, since this is a SVD scheme, the PMAs provide dual collections of orthonormal sets. 
\begin{theorem}
The principal moment functionals $\{ u_k({\bf x}) \}_{k=1}^r$ constitute an orthonormal set in $L^2({\bf R}^p, \mu)$.
\end{theorem}
\proof 

\begin{eqnarray}
\int_{{\bf R}^p} u_k({\bf x}) u_j({\bf x}) d\mu ({\bf x}) = \lambda_k^{-\frac{1}{2}} \lambda_j^{-\frac{1}{2}} \int_{{\bf R}^p} {\bf v}_j^T {\bf x} {\bf x}^T {\bf v}_k d \mu ({\bf x}) \nonumber \\
=  \lambda_k^{-\frac{1}{2}} \lambda_j^{-\frac{1}{2}} {\bf v}_j^T {\bf M}_2 (\mu)  {\bf v}_k = \delta_{k j} \; ,
\end{eqnarray}
where $\delta_{k j}$ is the Kronecker delta.
\endproof
We can complement these dual sets of PMAs to {\it orthogonal bases} of ${\bf R }^p$ ($\left\{ {\bf v}_k \right\}_{k=1}^p$) and $L^2({\bf R}^p, \mu )$ respectively. 

With $V$ denoting the orthogonal $p \times p$ matrix $V := [{\bf v}_1, \dots , {\bf v}_p ] $ we have that
\begin{equation}
V^T {\bf M}_2 (\mu) V = D
\end{equation}
where $D$ is the diagonal positive semidefinite matrix having the eigenvalues, {\it the principal moments}, of ${\bf M}_2 (\mu)$ on the diagonal, including the $p-r$ zeros representing the null space. In the PMA analysis framework, this spectral decomposition of ${\bf M}_2(\mu)$ and the accompanying dual sets of PMAs, allow us to define and explore optimal lower dimensional representations of the measure $\mu$ as well as associated representations of the original data, see https://principalmomentanalysis.github.io/ for an example of simplex based PMA.

We can also get exact estimates on how much of the total second moment of a measure that is captured in optimal lower dimensional representations. 

In fact, for any projection $\Pi \in  \mathcal{P}^s({\bf R}^p )$, from the linearity of the integral, we have that
\begin{equation}
\text{trace} ({\bf M}_2(\Pi_*(\mu))) = \text{trace} ( \Pi {\bf M}_2(\mu) \Pi ) 
\end{equation}
and so 
\begin{eqnarray}\label{eq:N=2 optimum}
\sup_{\Pi  \in \mathcal{P}^s({\bf R}^p)} \| {\bf M}_2(\Pi_*(\mu)) \|_{Tr} =  \sup_{\Pi  \in \mathcal{P}^s({\bf R}^p)} \text{trace} ( \Pi {\bf M}_2(\mu) \Pi ) = \nonumber \\
 \sup_{\Pi  \in \mathcal{P}^s({\bf R}^p)} \text{trace} ( V^T \Pi V V^T {\bf M}_2(\mu) V V^T  \Pi V ) = \nonumber \\
  \sup_{\Pi  \in \mathcal{P}^s({\bf R}^p)} \text{trace} ( \Pi  D  \Pi )  =  \sum_{k=1}^s \lambda_k \, ,
\end{eqnarray}
where the optimum is attained by a projection, $\Pi_s$, onto the eigenspace spanned by the $s$ "first" eigenvectors, PMAs,  of $M_2(\mu)$. This solution is unique if this eigenspace is uniquely defined, i.e. if the $s$ "first" eigenvectors are uniquely defined. This is the case unless the smallest concerned eigenvalue (possibly including $0$) happens to have higher geometric multiplicity than what is needed to define the projection. 

We conclude by noting that the second moment of a measure can be used, together with the total variation measure, to construct a natural pseudo-metric when measuring distances between measures and that, using the SVD scheme, we can give exact estimates of how much of the second moment we loose when we compare the original measure with optimal projections of the measure.

In fact, given two compactly supported Borel measures $\mu$ and $\nu$ on ${\bf R}^p$, let $|\mu - \nu |$ denote the {\it total variation measure} of $\mu - \nu $.
 We define the second moment semi-norm on compactly supported Borel measures to be the trace norm of  $ M_2(| \mu - \nu |)$, i.e.
\begin{equation}
\| \mu - \nu \|_2 := \| M_2(| \mu - \nu |) \|_{Tr} \, .
\end{equation}

Since, for any projection $\Pi$, we have that $\mu - \Pi_* \mu = (I-\Pi)_* \mu$ is a positive measure and
\begin{equation}
 \| {\bf M}_2( \mu) \|_{Tr} =  \| {\bf M}_2( \mu -\Pi_*(\mu)) \|_{Tr} + \| {\bf M}_2(\Pi_*(\mu)) \|_{Tr} \; ,
\end{equation}
we note that, with $\Pi_s$ denoting a rank $s$ solution to the optimization problem (\ref{eq:N=2 optimum}), we also get the exact error estimate
\begin{equation}
\inf_{\Pi  \in \mathcal{P}^s({\bf R}^p)} \| {\bf M}_2(| \mu -\Pi_*(\mu)|) \|_{Tr} = \| {\bf M}_2(| \mu -(\Pi_s)_*(\mu)|) \|_{Tr} 
=\sum_{k=s+1}^r \lambda_k \, .
\end{equation}


\begin{thebibliography}{999}
	\bibitem{pea} Pearson, K. (1901). "On Lines and Planes of Closest Fit to Systems of Points in Space". Philosophical Magazine. Series 6 Volume 2 (11): pp. 559–572
	\bibitem{ho1} Hotelling, H. (1933). "Analysis of a complex of statistical variables into principal components". Journal of Educational Psychology, 24, pp. 417–441, and 498–520.
	
\end{thebibliography}
\end{document}